\documentclass[%
 aip,
 jmp,%
 amsmath,amssymb,
 reprint,%
]{revtex4-1}

\usepackage{graphicx}
\usepackage{dcolumn}
\usepackage{bm}
\usepackage[mathlines]{lineno}
\begin{document}

\preprint{APS/123-QED}

\title{A Different Cell Size Approach to Fast Full-Waveform Inversion of Seismic Data}

\author{Amila Sudu Ambegedara}
 \email{nisanthaamila@gmail.com}
\affiliation{Department of Physical Sciences, Faculty of Applied Sciences, Rajarata University of Sri Lanka}%

\author{Indika Gayani Kumari Udagedara}
 \email{gayaniu@sci.pdn.ac.lk }
\affiliation{Department of Mathematics, Faculty of Science, University of Peradeniya, Sri Lanka}%


\begin{abstract}
\section*{Abstract}
Understanding the causes of sinkholes and determining the earth's subsurface properties will help Engineering Geologists in designing and constructing different kinds of structures. Also, determining of subsurface properties will increase possibilities of preventing expensive structural damages as well as a loss of life. Among the available health monitoring techniques, non-destructive methods play an important role. Full-waveform inversion together with the Gauss-Newton method, which we called as the regular method, able to determine the properties of the subsurface data from seismic data. However, one of the drawbacks of the Gauss-Newton method is a large memory requirement to store the Jacobian matrix. In this work, we use a different cell size approach to address the above issue. Results are validated for a synthetic model with an embedded air-filled void and compared with the regular method. \\~\\

\section*{Highlights}
\begin{itemize}
\item Full seismic waveform method based on Gauss-Newton method was used to detect embedded sinkholes in Earth's subsurface.
\item The difference cell size method is proposed to address the computational and memory requirements in Regular Full-wave inversion method.
\item Results are compared with regular full waveform inversion method
\item Less computational time is required for sinkhole detection with the proposed method.
\end{itemize}

\end{abstract}

                              
\maketitle


\section{\label{sec:level1}Introduction}

Determining of subsurface properties will help Engineering Geologists to prevent many Geo Hazards. Anomalies such as voids in soils cause significant structural damage. When a void weakens the support of the overlying earth, ground-surface depressions occur. Such a depression formed as a result of collapse is called a sinkhole. Engineering Geologists can use subsurface properties in designing human developments and constructing different types of structures in such away that minimizing future hazards.   

These collapses can result in significant property damage as well as a loss of life. For example, in 2013 in Florida, a man was swallowed by a sinkhole that opened beneath the bedroom of his house. This man's remains were never recovered. This sinkhole reopened up in 2015.  Also, repairing such damaged structures after a collapse is expensive. Thus, understanding the causes of sinkholes has the potential to prevent such expensive structural damage ahead of time. 

Identifications of sinkholes have been studied in the literature \cite{pazzi2018integrated, thomas1999evaluation,argentieri2015early} Methods based on Geo technical surface exploratory procedures such as cone penetration test and standard penetrations tests also have been used for evaluation of site characteristics\cite{nam2020preliminary, chang1999delineation, zini2015multidisciplinary}. Geologists have developed many testing methods for health monitoring in the geological sites \cite{sevil2017sinkhole}. Among them, non-destructive testing methods play an important role \cite{sevil2020characterizing}. There are many nondestructive testing methods available for sinkhole detection in a geological site.  Gravity methods \cite{wenjin1990effectiveness,gravity,mariita2007gravity}, electric resistivity methods \cite{van2002detection,resistivity}, and seismic methods\cite{nunziata2009s} are some exciting new methods in locating sinkholes. These methods have advantages and disadvantages in characterizing sinkholes. 

The full-waveform inversion (FWI) approach \cite{plessix2008introduction, virieux2009overview} is another approach that offers the potential to produce higher-resolution imaging of the subsurface by extracting information contained in the complete waveforms  \cite{tran2013sinkhole}. This approach can be determined the properties of the subsurface from seismic data (wavefield data) obtained at receivers, which are placed on the subsurface. 

In FWI, the model we consider is an initial guess based known properties of the subsurface. The model's results (wavefield data) are solved using wave equations assuming an elastic media producing ``predicted data". In our previous work \cite{ambegedara2021spatial}, a spatial mesh refinement method using cubic smoothing spline interpolation was proposed for forward modeling of FWI. 
Simultaneously, wavefield data is observed experimentally at the receivers, which are placed on the surface. Then the difference between the predicted data and the observed data is minimized to obtain properties of the subsurface. The model is updated iteratively until the residual is sufficiently small. 

In the FWI approach, the process of estimating wavefield by solving wave equations at known model parameters is known as the forward problem. If $M$ is the model space (or parameter space) and $D$ is the data (wave field) space, then the forward model $F: M \rightarrow D$ can be defined by
\begin{equation}
F\left(\textbf{m}\right) = \textbf{d}, 
\end{equation}
where $\textbf{m} \in M$ is the model parameters that represent the subsurface. For example, in the acoustic case, the model parameters are the P-wave velocities, S-wave velocities, density, and Lame coefficients defined at each cell of the numerical mesh used in the forward simulations.  $\textbf{d} \in D$ represents seismic responses of the surface recorded at the receivers. $F$ is the corresponding modeling operator, which is specified by the equation of motion and boundary conditions. 
Wavefield data obtained by forward simulation of wave equations and the observed seismic data. 

Finding $\textbf{m}$ by seeking the minimum of the residuals between the model responses obtained by simulation of wave equations and the observed seismic data is known as the inverse problem. The residual can be defined as
\begin{equation}
\Delta \textbf{d} = \textbf{d}_{\text{est}}\left(\textbf{m}\right) 
\end{equation}
where $\textbf{d}_{\text{est}}$ is the estimated data associated with the model parameters $m$ and $\textbf{d}_{\text{obs}}$ is observed data.


Recently, Ref.~\cite{tran2012site} developed an FWI technique that inverted body and surface waves in the case of real experimental data. This approach uses a Gauss-Newton technique to invert the full seismic wave-fields of near-surface velocity profiles by matching the observed and computed wave-forms in the time domain. Virtual sources and a reciprocity principle are used to calculate partial derivative wave-fields (gradient matrix) to reduce the computing time. 

The Gauss-Newton method consists of the computation of the Jacobian matrix, which records the partial derivatives of the seismic data. One of the drawbacks of the Gauss-Newton method is a large memory requirement to store the Jacobian matrix \cite{hu2011preconditioned, akcelik2002parallel,li2011compressed,tran2012site}.

In the literature, several techniques were used to reduce the memory usage for the Jacobian matrix \cite{hu2011preconditioned, akcelik2002parallel,li2011compressed}. For example, Ref. \cite{hu2011preconditioned} used a non-linear conjugate gradient method for seismic wave inversion as it does not require the inversion of the dense Hessian matrix. However, the convergence rate of the results may be slow with the conjugate gradient method and not efficient for the problems with more parameters. Ref. \cite{akcelik2002parallel} proposed a Gauss-Newton-Krylov based method, which is a matrix-free implementation of the Gauss-Newton method for full-wave inversion problems. The authors in Ref. \cite{akcelik2002parallel} showed that this approach is well suited for
a nonlinear and ill-conditioned problem such as inverse wave propagation. A compressed implicit Jacobian scheme for 3D electromagnetic data inversion was proposed in Ref. \cite{li2011compressed}.  A significant reduction in memory usage for the Jacobian matrix is obtained with the implicit Jacobian scheme for reconstructing electromagnetic data.

In this work, we introduce a different cell size based technique as an option for the Jacobian matrix storage. The goal is to address the computational efficiency and the memory requirements of the developed method in Ref.~\cite{tran2012site}. The difference cell size based technique is applied to a synthetic model and compared with the Gauss-Newton inversion regular method, which is developed in Ref.~\cite{tran2012site}. 

The paper is organized as follows. The wave propagation equations in elastic media is presented in Section
 \ref{der}. We present the Full-wave inversion method, which was introduced in Ref.~\cite{tran2012site} in Section \ref{FWI}. Section \ref{cell_method} presents the different cell size approach. The comparison of the methods and results for a synthetic model are presented in Section \ref{results}. 

\section{Wave Propagation in Elastic Media}\label{der}
Equations of wave propagation in elastic media are derived by using Newton's Second Law of Motion and Hooke's Law \cite{wave_eq}. These equations can be derived by considering the total force applied to a volume element of an elastic media. 
We can express the equations of wave propagation as
\begin{align}
\rho\frac{ \partial v_i}{\partial t} &= f_i + \frac{ \sigma_{ij}}{\partial x_j} \label{prop_1} \\
\frac{\partial \sigma_{ij}}{\partial t} &= \lambda \frac{\partial \theta}{\partial t} \delta_{ij} - 2 \mu \frac{\partial \epsilon_{ij}}{\partial t} \label{prop_2} \\
\frac{\partial \epsilon_{ij}}{\partial t} &= \frac{1}{2} \left( \frac{ \partial v_i}{ \partial x_j} + \frac{\partial v_j }{\partial x_i} \right) \label{prop_3},
\end{align}
where $v$ is the particle velocity, $\sigma_{ij}$ is the shear stress, $\lambda$ and $\mu$ ate two elasticity coefficients,which are called the Lame parameters, $\theta$ is the stress tensor,and $\epsilon_{ij}$ is strain tensor.

Seismic waves are elastic waves. The two independent parameters in elastic tensor can be expressed in terms of elastic moduli.  If $\kappa$ represents the bulk modulus of the material, then  
\begin{equation}
 \lambda = \kappa - \frac{2 \mu}{3}
\end{equation}
and the Young's moduli
\begin{equation}
E = \frac{\left(3 \lambda + 2 \mu \right) \mu}{\left( \lambda + \mu \right)}.
\end{equation}
The wave propagation velocity depends on the elasticity and the density of the medium. The P-wave and S-wave velocities are
\begin{equation}
V_p = \sqrt{\frac{\kappa + 4 \mu /3}{\rho}} = \sqrt{\frac{ \lambda + 2 \mu }{\rho}}
\end{equation}
and
\begin{equation}
V_s = \sqrt{\frac{\mu}{\rho}},
\end{equation}
where $V_s$ and $V_p$ are P-wave and S-wave velocities of the medium. Moreover, the shear moduli G is defined as 
\begin{equation}
G = \mu = \rho V_s^2. 
\end{equation}

\section{Full-Waveform Inversion} \label{FWI}

The FWI technique consists of two stages. The first stage includes forward modeling to generate synthetic
wave-fields and the second stage includes the model updating by considering when the residual between predicted and measured surface velocities are negligible.
In this thesis, we consider wave equations in 2-D cartesian coordinates.
\subsection{Forward Problem}
Forward modeling seeks the solutions of the 2-D elastic wave equations. We simulate wave propagation by solving 2-D elastic wave equations numerically. The governing equations for 2-D elastic wave propagation can be obtained using Equations \ref{prop_1}-\ref{prop_3}. 

Let $ \sigma_{xx}$,  $\sigma_{zz}$, and $\sigma_{xz}$ be the components of stress tensor and $u$, $v$ be the particle velocity components. The spatial directions in the 2D plane are $x$ and $z$.

Then the equations governing particle velocity in 2-D are
\begin{align}
 \frac {\partial u}{\partial t} = \frac{1}{\rho}\left( \frac {\partial \sigma_{xx}}{\partial x} + \frac {\partial \sigma_{xz}}{\partial z}\right)=f_1 \left(\rho \right) \label{f1}\\
 \frac {\partial v}{\partial t} = \frac{1}{\rho} \left( \frac {\partial \sigma_{xz}}{\partial x} + \frac {\partial \sigma_{zz}}{\partial z} \right) = f_2 \left(\rho \right) \label{f2}
\end{align}
and the equations governing stress-strain tensor are 
\begin{align}
 \frac {\partial \sigma_{xx}}{\partial t} =& \left(\lambda +2 \mu \right) \frac {\partial u}{\partial x} +  \lambda \frac {\partial v}{\partial z} = f_3 \left(\lambda, \mu \right) \label{f3}\\
  \frac {\partial \sigma_{zz}}{\partial t} =& \lambda \frac {\partial u}{\partial x} + \left(\lambda +2 \mu \right) \frac {\partial v}{\partial z}= f_4 \left(\lambda, \mu \right)\label{f4} \\
 \frac {\partial \sigma_{xz}}{\partial t} =& \mu  \left( \frac {\partial v}{\partial x} +  \frac {\partial u}{\partial z}\right)=f_5 \left(\mu \right) \label{f5}
\end{align}
Here $\rho \left(x, z\right)$  is the mass density, $\mu  \left(x, z\right)$, and $\lambda  \left(x, z\right)$ are the Lame's coefficients of the material. 
The equations \ref{f1}-\ref{f5} can be written as 
\begin{equation}\label{fwd}
 F\left( \rho \left(x, z\right), \mu  \left(x, z\right), \lambda  \left(x, z\right) \right) = \textbf{d}.
\end{equation}
Equations \ref{f1}-\ref{f5} are the forward equations of the FWI method. We can express the forward equations in the form of $A \textbf{x} = \textbf{b}$, where 
\begin{widetext}
\begin{equation}
A = \begin{bmatrix}
\frac {\partial \sigma_{xx}}{\partial x} + \frac {\partial \sigma_{xz}}{\partial z} & 0 & 0 \\
\frac {\partial \sigma_{xz}}{\partial x} + \frac {\partial \sigma_{zz}}{\partial z} & 0 & 0 \\
0 &  \frac {\partial u}{\partial x} + \frac {\partial v}{\partial z} & 2 \frac {\partial u}{\partial x} \\
0 &  \frac {\partial u}{\partial x} + \frac {\partial v}{\partial z} & 2 \frac {\partial v}{\partial x} \\
0 & 0 &  \frac {\partial v}{\partial x} +  \frac {\partial u}{\partial z} \\
\end{bmatrix}, \textbf{x} = \begin{bmatrix} 
\frac{1}{\rho} \\
\lambda \\
\mu
\end{bmatrix}, \text{and} ~ \textbf{b} = \begin{bmatrix}  
\frac {\partial u}{\partial t}  \\
 \frac {\partial v}{\partial t} \\
 \frac {\partial \sigma_{xx}}{\partial t} \\
 \frac {\partial \sigma_{zz}}{\partial t} \\
 \frac {\partial \sigma_{xz}}{\partial t}
\end{bmatrix}
\end{equation}
\end{widetext}

To solve the above forward equations numerically, specific boundary conditions are needed. 
We impose three boundary conditions: the free surface boundary condition on the top of the domain, the absorbing boundary condition on the right side of the domain and bottom of the domain, and the symmetric boundary condition on the left-hand side of the domain. 

\subsubsection{Free Surface Boundary Conditions}

The measurements of the wavefield are generally collected along the earth's subsurface. Therefore, we impose the free surface boundary condition on the top of the domain by setting the vertical stress components are as zero. 

\begin{equation}\label{bd1}
\begin{cases}
 \sigma_{xz} = 0 \\
 \sigma_{zz} = 0.
\end{cases}
\end{equation}

\subsubsection{Absorbing Boundary Conditions}
Numerical methods are solved for a region of space by imposing artificial boundaries. Therefore, to avoid the reflections from the boundaries, absorbing boundary conditions should be applied on the right-hand side and the bottom of the domain. Thus the absorbing condition at the bottom of the domain is 

\begin{equation}\label{bd3}
\begin{cases}
\frac {\partial u}{\partial t} + V_s \frac {\partial u}{\partial z}  = 0 \\
\frac {\partial v}{\partial t} + V_p \frac {\partial v}{\partial z}  = 0
\end{cases}
\end{equation}
and at the right-hand side of the domain
\begin{equation}\label{bd4}
\begin{cases}
\frac {\partial u}{\partial t} + V_s \frac {\partial u}{\partial x}  = 0 \\
\frac {\partial v}{\partial t} + V_p \frac {\partial v}{\partial x}  = 0,
\end{cases}
\end{equation}
where $V_s$ and $V_p$ are sheer and pressure wave velocities, respectively. 

\subsubsection{Symmetric Condition}
To save computational time, we imposed a symmetric condition along the load line. Thus at the left-hand side of the domain we set

\begin{equation}\label{bd2}
\begin{cases}
 \sigma_{xz} = 0 \\
 u = 0.
\end{cases}
\end{equation}

To solve the forward equations, one can use numerical approaches such as finite difference method, finite element method, and Fourier/spectral method. 
Ref.~\citep{tran2012site} used a classic velocity-stress staggered-grid finite-difference solution of the
2-D elastic wave equations in the time domain (Virieux, 1986)  with the absorbing boundary conditions (Clayton
and Engquist, 1977).  In that approach, a direct discretization of the equations \ref{f1}-\ref{f5}, both in time and in space is considered. We follow the same approach for solving forward equations.

\subsection{ A Classic Finite Difference Scheme}\label{finite}

To  solve Eqs.~\ref{f1}-\ref{f5} with the above boundary conditions \ref{bd1} - \ref{bd4}, the derivatives are discretized using central finite differences. 

\begin{figure}
\centering
\includegraphics[width=1.5\linewidth]{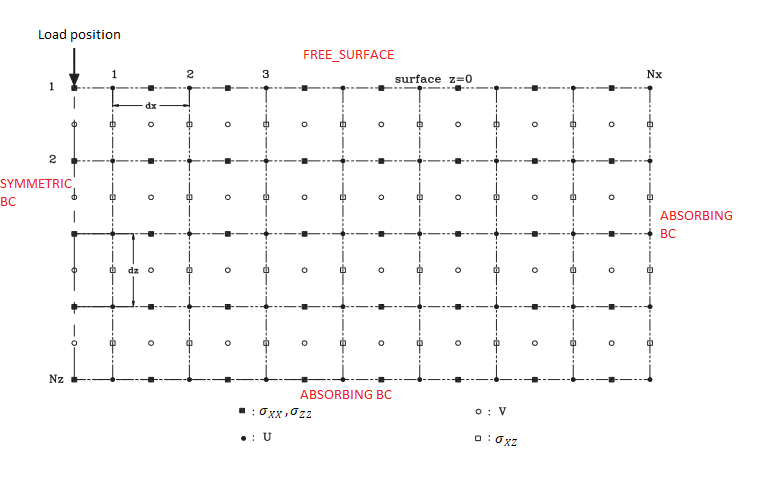}
\caption{The discretization of the domain  }
\label{dom_fig}
\end{figure}

In our problem, for a field variable $f$, the temporal discretization is
\begin{align}
D_t \left[ f \right]_{i,j}^{k} = \frac{f_{i,j}^{k+1/2} - f_{i,j}^{k-1/2}}{\delta t} = \frac{\partial f}{\partial t} \arrowvert_{i,j}^{k} + O(\delta^2) 
\end{align}
and the spatial discretizations are
\begin{align}
D_x \left[ f \right]_{i,j}^{k} =& \frac{f_{i+1/2,j}^{k} - f_{i-1/2,j}^{k}}{h_1} = \frac{\partial f}{\partial x} \arrowvert_{i,j}^{k} + O(h_1^2) \\
D_z \left[ f \right]_{i,j}^{k} =& \frac{f_{i,j+1/2}^{k} - f_{i,j-1/2}^{k}}{h_3} = \frac{\partial f}{\partial z} \arrowvert_{i,j}^{k} + O(h_3^2),
\end{align}
where $\textbf{O} \left(\cdot\right)$ is the truncation error. Here $i, j$, and $k$ represent the indices used in the discretization for the directions $x, y$ and time. The domain is discretized in the $x, y$ and time directions as shown in Fig.~\ref{dom_fig}. $h_1, h_3$, and $\delta t$ are the grid steps for $x$, $z$ and time directions, respectively. $f$ can take $u, v, \sigma_{xx}, \sigma_{zz},\sigma_{xz}$.
For example, the derivative terms $\frac{\partial u}{\partial t}$, $\frac{\partial \sigma_{xx}}{\partial x}$, and $\frac{\partial \sigma_{xz}}{\partial z}$ in  Eq.~\ref{f1} can be approximated as
\begin{align} 
\frac{\partial u}{\partial t} &= \frac{u_{i,j}^{k+ 1/2} - u_{i,j}^{k- 1/2}}{2 \delta t} \label{d1} \\
\frac{\partial \sigma_{xx}}{\partial x} &= \frac{{\sigma_{xx}}^{k}_{i+ 1/2,j} - {\sigma_{xx}}^{k}_{i- 1/2,j}}{2 h_1} \label{d2} \\
\frac{\partial \sigma_{xz}}{\partial z} &= \frac{{\sigma_{xz}}^{k}_{i,j+1/2}- {\sigma_{xz}}^{k}_{i,j-1/2}}{2 h_3} \label{d3} 
\end{align}
Then, Eq.~\ref{f1} can be approximated using Eqs. \ref{d1}, \ref{d2}, and \ref{d3} as,
\begin{widetext}
\begin{align}
 \frac{u_{i,j}^{k+ 1/2} - u_{i,j}^{k- 1/2}}{2 \delta t} = \frac{1}{\rho} \left( \left( \frac{{\sigma_{xx}}_{i+ 1/2,j}^{k} - {\sigma_{xx}}_{i-1/2,j}^{k}}{2 h_1} \right) \left(  \frac{{\sigma_{xz}}_{i,j+1/2}^{k} - {\sigma_{xz}}_{i,j-1/2}^{k}}{2 h_3}\right)  \right)
\end{align}
\end{widetext}

Equations~\ref{num1} - \ref{num5} are the second order accuracy numerical scheme after discretizing the system of differential equations (Virieux, 1986). The velocity field $\left(U, V\right) = \left(u, v\right)$ at time $\left(k+\frac{1}{2}\right) \delta t$ and the stress-tensor field $\left(T_{xx}, T_{zz}, T_{xz} \right) = \left( \sigma_{xx}, \sigma_{zz}, \sigma_{xz} \right) $ at time $\left(k+1\right) \delta t$ are explicitly calculated with the numerical scheme. 
\begin{widetext}
\begin{align} 
U_{i,j}^{k +1/2} &= U_{i,j}^{k-{1/2} }+ B_{i,j} \frac{\delta t}{h_1} \left( Txx_{i + {1/2},j}^k -  Txx_{i - {1/2},j}^k \right)  + B_{i,j} \frac{\delta t}{h_3} \left( Txz_{i, j +{1/2}}^k -  Txz_{i, j -{1/2}}^k \right) \label{num1} \\
V_{i+{1/2},j+{1/2}}^{k +1/2} &=V_{i+{1/2},j+{1/2}}^{k-{1/2} }+ B_{i+{1/2},j+{1/2}} \frac{\delta t}{h_1} \left( Txz_{i + 1,j+{1/2}}^k -  Txz_{i ,j+{1/2}}^k \right) \nonumber \\
&+ B_{i+{1/2},j+{1/2}} \frac{\delta t}{h_3} \left( Tzz_{i+{1/2}, j +1}^k -  Txz_{i+{1/2}, j}^k \right) \label{num2}\\
Txx_{i+{1/2},j}^{k +1} &=Txx_{i+{1/2},j}^{k }+ \left(L+2M\right)_{i+{1/2},j} \frac{\delta t}{h_1} \left( U_{i + 1,j}^{k+{1/2}} -  U_{i ,j}^{k+{1/2}} \right) \nonumber \\
&+ L_{i+{1/2},j} \frac{\delta t}{h_3} \left( V_{i+{1/2}, j +{1/2}}^{k+{1/2}} -  U_{i+{1/2}, j -{1/2}}^{k+{1/2}} \right) \label{num3}\\
Tzz_{i+{1/2},j}^{k +1} &=Tzz_{i+{1/2},j}^{k}+ \left(L + 2M \right)_{i+{1/2},j} \frac{\delta t}{h_1} \left( V_{i + {1/2},j+{1/2}}^{k+{1/2}} -  V_{i+{1/2} ,j-{1/2}}^{k+{1/2}} \right) \nonumber \\
&+ L_{i+{1/2},j} \frac{\delta t}{h_3} \left( U_{i+1, j}^{k+{1/2}} - U_{i, j}^{k+{1/2}} \right) \label{num4}\\
Txz_{i,j+{1/2}}^{k +1} &=Txz_{i,j+{1/2}}^{k }+ M_{i,j+{1/2}} \frac{\delta t}{h_3} \left( U_{i ,j+1}^{k+{1/2}} -  U_{i ,j}^{k+{1/2}} \right) \nonumber \\
&+ M_{i,j+{1/2}} \frac{\delta t}{h_1} \left(V_{i+{1/2}, j +{1/2}}^{k+{1/2}} -  V_{i-{1/2}, j +{1/2}}^{k+{1/2}} \right) \label{num5}
\end{align}
\end{widetext}
 Here, $M$ and $L$ represent the Lame coefficients $ \left(\mu, \lambda\right)$ and 
\begin{equation}
B = \frac{1}{\rho}
\end{equation}
as shown in  Fig.~\ref{dom_fig}. 

Moreover, the initial condition at time $t=0$ is set such that the stress and velocity are zero everywhere in the domain. The medium is perturbed by changing vertical stress $\sigma_{zz}$ at the source using 
\begin{equation}\label{ini}
R\left(t\right) = \left[1- 2 \pi ^2 f_c^2 \left(t -t_0\right)^2 \right] \exp \left[ - \pi^2 f_c^2 \left(t -t_0\right)^2 \right],
\end{equation}
where $f_c$ is the center of the frequency band and $t_0$ is the time shift. 

\subsubsection{Stability Criterion}

Numerical schemes are associated with numerical errors due to the approximation of the derivatives in the partial differential scheme. It is important to obtain a stable wave propagation solution from the finite difference scheme. With some numerical schemes, the errors made at one-time step grow as the computations proceed. Such a numerical scheme is said to be unstable so the results blow up. If the errors decay with time as the computations proceed, we say a finite difference scheme is stable. In that case, the numerical solutions are bounded. 

To obtain a bounded solution from the finite difference scheme, we obtain $\delta t$ from the stability criterion (Virieux, 1986) given by
\begin{equation}
\delta t  \leq \frac{1}{V_{\text{max}} \sqrt{\frac{1}{h_1^2} + \frac{1}{h_3^2}}}.
\end{equation}
Here $V_{\text{max}}$ is the maximum P-wave velocity in the media. 

Inputs for the forward problem are the model parameters such as density, Lames's moduli, P-wave velocity, and S-wave velocity. Then the particle velocities and stresses (outputs) are calculated by implementing the numerical scheme (Eqs.~\ref{num1} - \ref{num5}) in MATLAB. 

\subsection{Inverse Problem}\label{inv}
The FWI is the problem of finding the parametrization of the subsurface using seismic wave field. Thus the goal of inversion is to estimate a discrete parametrization of the subsurface by minimizing the residual between the observed seismic data and the numerically predicted seismic data. 
If seismic waves are generated from $NS$ sources (one shot at a time) and are recorded by $NR$ receivers, then
 the residual for all shots and receivers can be defined as 
\begin{equation}
\Delta \textbf{d}_{ij} = F_{ij} \left(\textbf{m}\right) - \textbf{d}_{ij},
\end{equation}
where $\textbf{d}_{i,j}$ and $F_{i,j}\left(\textbf{m}\right)$ are the observed data and the estimated data
associated with the model parameters $\textbf{m}$, and indices $i$ and $j$ denote the $i^{th}$
shot and $j^{th}$ receiver, respectively. In this problem, the model parameters are density, $\rho \left(x, z\right)$, or one of Lame's moduli, $\lambda \left(x,z\right)$ and $\mu\left(x, z\right)$ or $V_s \left(x, z\right)$ and $V_p\left(x, z\right)$. However, due to the relationship between elastic moduli and wave velociies only three model parameters are enough to characterize the subsurface.

This problem can be modeled as a least squares problem. The objective function of the inverse problem can be expressed as minimizing a least square error $E_d \left(m\right)$. For $\Delta \textbf{d} : M \rightarrow D$ the problem is
\begin{equation}
\text{argmin}_{\textbf{m}} E_d \left(\textbf{m} \right) = \frac{1}{2} \Delta \textbf{d}^T \Delta \textbf{d} = \frac{1}{2} || \Delta \textbf{d} ||_2^2,
\end{equation}
where $E_d \left(\textbf{m}\right)$ is called the misfit function \cite{sheen2006time}. Here $\Delta \textbf{d} = \{ \Delta \textbf{d}_{i,j}, i=1,..., NS, j=1,...,NR\}$. $\Delta \textbf{d}$ is a
column vector, which is the combination of residuals $\Delta \textbf{d}_{i,j}$ for all
shots and receivers. Optimization problems of this form are called nonlinear least-squares problem and our target here is to find model parameters $\textbf{m}^*$ that minimizes 
$E_d \left(\textbf{m}\right)$. 

Model updating methods such as Gradient descent method, Newton method, and the Gauss-Newton method can be used to solve the above optimization problem. 

\subsubsection{Gradient Descent  Method}
The gradient method solve the non linear least square problem with search directions defined by the gradient of the function $E_d\left(\textbf{m}\right)$. $E_d \left(\textbf{m}\right)$ decreases in the negative direction of the gradient of $E_d\left(\textbf{m}\right)$, $-\nabla E_d \left(\textbf{m}\right)$. For iterations  $\text{n} \geq 0$,
\begin{equation}
\textbf{m}^{n+1} = \textbf{m}^{n} - \alpha \nabla E_d \left(\textbf{m}\right)
\end{equation}
iterates to find the minimum number $\textbf{m}^*$. Here $\alpha$ is the step size. The gradient contains the first partial derivatives of 
$E_d \left(\textbf{m}\right)$ with respect to the model parameters $\textbf{m}$. 
\begin{equation}
\nabla E_d \left(\textbf{m}\right) = \left[ \frac{ \partial F_i \left(\textbf{m}\right)}{\partial m_j} \right]^T \Delta d,
\end{equation}
where $i=1,...,NS\times NR$ and $j=1,..., M$. Here $\left[ \frac{ \partial F_i \left(\textbf{m}\right)}{\partial m_j} \right]$ is defined as the Jacobian matrix $J$. 

\subsubsection{Gauss-Newton Method}

The Newton method is based on the model update with the second order partial derivatives of the function $E_d \left(\textbf{m}\right)$. 
The Hessian matrix, 
\begin{equation}
H = -\nabla^2  E_d \left(\textbf{m}\right) =  \frac{\partial}{\partial m_p} \left[J^T \Delta \textbf{d} \right]
\end{equation}
records the second order derivatives. The Newton method is given by
\begin{equation}
\textbf{m}^{n+1} = \textbf{m}^n - H^{-1} J^T \Delta \textbf{d}. 
\end{equation}
However, the calculation of the Hessian matrix is difficult. Therefore, the Newton method has not been often used in geophysical inverse problems. 
The Hessian matrix can be written as 
\begin{equation}\label{he}
H = J^T J + \frac{\partial J^T}{\partial m_p} \Delta \textbf{d}.
\end{equation}
By considering the negligibility of the second term of the Eq. \ref{he}, the Hessian matrix can be approximated as
\begin{equation}
H_a =  J^T J
\end{equation}
and the Gauss-Newton formula  \cite{pratt1998gauss}  with the approximate Hessian matrix is 
\begin{equation}
\textbf{m}^{n+1} = \textbf{m}^n - H_a^{-1} J^T \Delta \textbf{d}. 
\end{equation}

The Gauss-Newton method is effective for solving non-linear problems and guarantees faster convergence rates than the gradient method. With good initial guesses, the Gauss-Newton method converges nearly quadratically. But theoretically, the Gauss-Newton method converges linearly. However, when the Jacobian is ill-conditioned or singular, the search direction becomes very large and the Gauss-Newton method is not globally convergent. Thus to solve the original problem, regularization of the original problem can be used. 

The regularized misfit function $E$ used in geophysical inversion\cite{sheen2006time} is defined as 
\begin{equation}
E \left( \textbf{m} \right) = E_d \left(\textbf{m}\right) + \lambda E_m \left(\textbf{m}\right),
\end{equation}
where $\lambda$ is the regularization parameter that controls the relative importance of the $E_m \left(\textbf{m}\right)$, where $E_m \left(\textbf{m}\right)$ is the model objective function that contains a priori information of the model. $E_m \left(\textbf{m}\right)$ can be written as 
\begin{equation}
E_m \left(\textbf{m}\right) = \frac{1}{2} || L \Delta \textbf{m} ||^2, 
\end{equation}
where $L$ is discrete linear operator \cite{sheen2006time}. Then the regularized Gauss-Newton formula with step size $\alpha$ is
\begin{equation}\label{reg_g1}
\textbf{m}^{n+1} =  \textbf{m}^{n} - \alpha^n \left[J^T J + \lambda L^T L \right]^{-1} J^T \Delta \textbf{d}.
\end{equation}
When $L=I$, Eq. \ref{reg_g1} represents the damped least-squares method \cite{sheen2006time}. $L$ can be used as a discrete 2-D Laplacian operator\cite{sasaki1989two}, which is defined as
\begin{widetext}
\begin{equation}
L_i \Delta \textbf{m} \approx P_i \Delta \textbf{m} = \left( \Delta m_i \right) ^E +  \left( \Delta m_i \right) ^W+ \left( \Delta m_i \right) ^N +\left( \Delta m_i \right) ^S - 4 \left( \Delta m_i \right),
\end{equation}
\end{widetext}
where $E, W, N,$ and $S$ are the four neighbors of the $i^{th}$ model parameter and $P_i$ is the $i^{th}$ row of the Laplacian matrix whose elements are either 1, -4, or 0. Ref. \cite{sheen2006time} used both model objective functions from damped least-squares method and discrete 2-D Laplacian operator in the regularized problem. The regularized Gauss-Newton formula \cite{sheen2006time} to geophysical inversion can be written as
\begin{equation}\label{update}
\textbf{m}^{n+1} =  \textbf{m}^{n} - \alpha^n \left[J^T J + \lambda_1 P^T P + \lambda_2 I^T I \right]^{-1} J^T \Delta \textbf{d}.
\end{equation}
The step length $\alpha^n$ is determined by 
\begin{equation}
\alpha^n \approx \frac{\left[ J^T g^n \right]^T \left[ F\left(\textbf{m}^n\right) -\textbf{d} \right]}{\left[ J^T g^n\right]^T \left[ J^t g^n \right]},
\end{equation}
where $g^n = \left[J^T J + \lambda_1 P^T P + \lambda_2 I^T I\right]^{-1} J^T \left[F\left(\textbf{m}^n\right) - \textbf{d}\right]$.

For model updating, Ref.~\cite{tran2012site} uses  Eq. \ref{update} in the FWI. Following modifications to the residual and the Jacobian matrix are also used in Ref.~\cite{tran2012site}. 

\begin{enumerate}
\item The residual $\Delta d_{ij}$ is modified to avoid the influence of the source on the estimation during inversion. For that modification cross-convolution of wave-fields is used.  The symbol * denotes the convolution. Let the model $\textbf{m}$ includes all unknowns (S-wave and P-wave velocities of cells). For each shot gathering, the estimated
wave-fields are convolved with a reference trace from the
observed wave-field, and the observed wave-fields are convolved
with a reference trace from the estimated wave-field. Thus the modified residual between estimated and observed data for the
$i^{th}$ shot and $ j^{th}$ receiver is
\begin{equation}
\Delta d_{ij} = F_{i,j} \left(\textbf{m}\right) * d_{i,k}- d_{ij}*F_{i,k} \left(\textbf{m}\right), 
\end{equation}
where $d_{i,j}$ and $F_{i,j} \left(\textbf{m}\right)$ are the observed data and the estimated
data associated with the model $\textbf{m}$. 
$F_{i,k}\left(m\right)$ and $d_{i,k}$ are the reference traces from the estimated and observed data, respectively, at the $k^{th}$ receiver position.

\item The Jacobian matrix $J$ is obtained by taking the partial derivatives of seismograms with respect to parameters of model $\textbf{m}$ and convolving with the reference traces and defined by
\begin{equation}\label{Jac}
J_{ij,M} = \frac{\partial F_{i,j} \left(\textbf{m}\right) }{\partial m_p} * d_{i,k} - d_{i,j} * \frac{\partial F_{i,k} \left(\textbf{m}\right)}{\partial m_p},
\end{equation}
for $i=1,..., NS, j=1,...,NR$, and $p =1,...,M$. 
\end{enumerate}
 



Here $\lambda_1$ and $\lambda_2$ are between 0 and $\infty$. For this study, $\lambda_1= 0.05$ and $\lambda_2 = 0.0005$ are chosen as appropriate values.

\section{Different Cell Size Method to Store Jacobian} \label{cell_method}

To solve the inverse problem introduced in Section \ref{inv}, the regularized Gauss-Newton formula, which is defined in Eq.~\ref{update}, can be used. 
The term $J^T J$ is defined as approximation to the Hessian matrix $H$. The major drawback of the Gauss-Newton method is memory and computational requirements of the Hessian matrix approximation. 

The size of the Jacobian matrix is equal to the number of measured data points at receivers for each shot ($NR \times NS$) times the number of parameters (number of cells in the domain). In 3-D problems the size of the data set and the number of cells in the domain is usually large. Hence, the storage of the Jacobian matrix requires an adequate amount of storage. For large scale problems, as the size of the Jacobian matrix increases, expenses to calculate $J^T J$ and invert $J^T J + \lambda_1 P^T P + \lambda_2 I^T I$ in Eq.~\ref{update} also increase.  As we discussed in the itroduction, there are several approaches to manage the storage and computational requirements of the inversion. Other than those techniques, 
Ref.~\cite{sheen2006time} suggested a way to calculate $ H_a$ matrix without fully storing the Jacobian matrix. In their approach, the Jacobian matrix was divided into sub-matrices at the receivers and $H_a$ matrix was calculated as follows:
\begin{equation}
H_a = J^T J =
\begin{bmatrix}
    J_1^T     &  J_2^T & \dots & J_{NR}^T 
\end{bmatrix}
\begin{bmatrix}
    J_1^T   \\  
J_2^T  \\
\dots \\
J_{NR}^T  \\
\end{bmatrix}=\sum_j^{NR} J_{j}^{T} J_{j}
\end{equation}
In this way, the $H_a$ can be calculated by summing up the sub-matrices. The above technique can be implemented in MatLab with a loop that goes through the number of receivers. Therefore, the full Jacobian matrix does not have to be stored. Ref.~\cite{tran2012site} used the same technique for calculating $H_a$ matrix. In the implementation, Ref.~\cite{tran2012site} did all the above manipulation with arrays rather than the matrices. Before calculating the $H_a$ matrix, Ref.~\cite{tran2012site} converted the Jocobian sub matrices to an array. 

In this work, we introduce an approach called ``different cell size method" in addition to the computational  techniques used in Ref.~\cite{sheen2006time} and  Ref.~\cite{tran2012site}. One special observation on the Jacobian matrix is that the partial derivatives values of the seismograms corresponds to the bottom cells in the domain are smaller when compared with those values at the top cells. By considering that fact, we decompose the spatial domain into three zones: zone 1, zone 2, and zone 3. One can rather choose more zones according to the size of the domain. Here the Jacobian matrix has the smallest values at the cells corresponding to the zone 3. Then zone 2 and 3 are discretized again with a bigger step size in the $x$ direction and $z$ direction. The discretization ratio for zone 3 is larger than that of zone 2 and the discretization ratio for zone 2 is larger than that of zone 1.  For example, if $dx$ and $dz$ are step size for the regular domain in the $x$ and $z$ direction, then the discretization ratios for zone 2 can be $2 dx$ and $2 dz$. The discretization ratios for zone 3 can be $3 dx$ and $3 dz$. According to that, one cell in the zone 2 is created by combining 4 smaller cells and one cell in the zone 3  is created by combining 9 smaller cells in the regular domain

\begin{figure}
\centering
\includegraphics[height=0.7\linewidth]{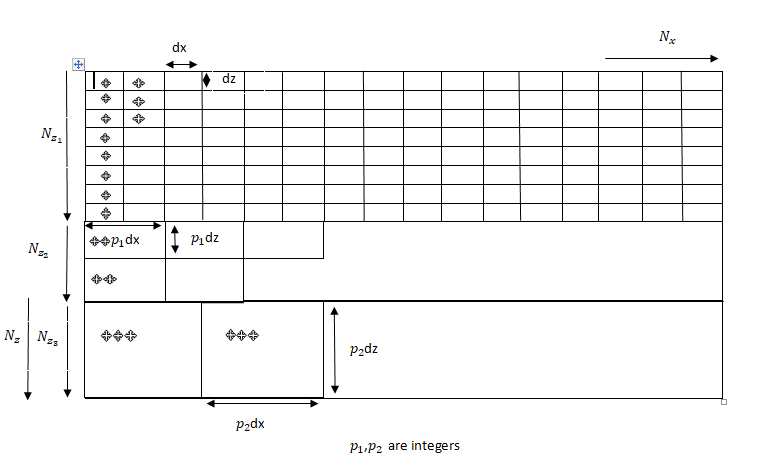}(a)
\includegraphics[height=0.8\linewidth]{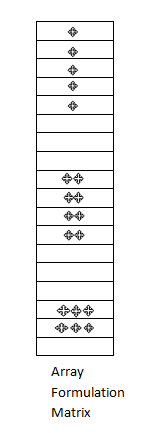}(b)
\caption{(a) Decomposition of the initial domain. The cell size of zone 1, zone 2, zone 3 are $dx \times dz$, $ p_1 dx \times p_2 dz$, and $p_3 dx \times p_3 dz$, respectively. Here $p_1, p_2$ are integers. (b) The array formulation for the values in the cells in the domain.}
\label{dom_cell}
\end{figure}

 The values of the Jacobian at the bigger cells in  zone 2 and zone 3 are re-evaluated by taking the sum of the values at the smaller cells. Then the values at the cells for the three zones  are stored in matrices and converted to a single array. Figure \ref{dom_cell}  illustrates the procedure of the new discretization and converting values from matrix to an array for the three zones.
Notice that the length of the obtained array for the different cell size method is shorter than the length of the array obtained with the regular cell method, which uses Jacobian matrix without combining cells, for the initial domain. For the regular cell size method, the number of parameters is equal to the number of cells in the domain. For the proposed different cell size method, the number of cells is less than that of the regular cell size method as the bottom cells in the bottom zones are combined. Thus the size of the Hessian approximation matrix is smaller than that of the regular cell size method. Due to that, the Hessian approximation matrix can be calculated faster and less Jacobian storage is required. Thus the different cell size method is computationally inexpensive compared with the regular cell size method. Once we calculate the Hessian matrix, the Gauss-Newton update was used to find the shear wave and pressure wave velocities. Then the velocities in the bigger cells (combined cells) in zone 2 and zone 3 are converted back to smaller cells.  The  velocities at the smaller cells are calculated by taking the average of the bigger cells. 

\section{Numerical Results} \label{results}

In this section we investigate the capacity of the FWI in detection of embedded voids. The FWI technique with the different cell size method is applied to a synthetic model. Results are compared with the regular cell size method, which was used in Ref.~ \cite{tran2012site}. 

\subsection{A Synthetic Model with an Embedded Void} 

We consider a synthetic model of the earth for the investigation. The velocity profiles of the earth, i.e., S-wave and P-wave velocities of cells, are assumed to be known a priori. In the test configuration, the locations of a set of sources and receivers are also known. Using a known velocity structure, surface waveform data are
calculated. These waveform data are then used as the  input to the inversion program. If the waveforms were obtained from a field
test, then the velocity structures can be extracted from the inversion of the surface waveform data. Theoretically, the extracted velocity profile
should be the same as the velocity profile assumed at the start.

\begin{figure}
\centering
\includegraphics[width=1\linewidth]{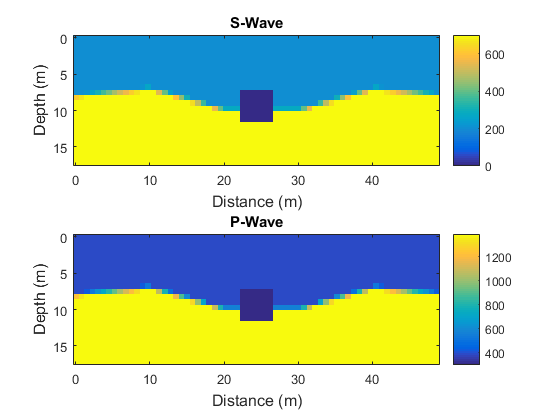}
\caption{Velocity profiles of the true model}
\label{true_fig}
\end{figure}

We consider a synthetic model, which consists of two layers with an embedded air-filled void. 
The S-wave velocities $V_s$ of the materials are 200 $\text{m/s}$ for the soil layers and 700 $\text{m/s}$ for limestone. 
The P-wave velocity
is generated from the S-wave velocity $V_p$ using 
\begin{equation}\label{pwave}
V_p=\sqrt {\left(2 \left(1-\nu\right)/\left(1-2 \nu \right) \right)} V_s
\end{equation}
for the entire domain. Here $\nu$ is 0.33.  The void is encoded by setting the S-wave velocity in some computational cells to zero and P-wave velocity of those cells to 300 $\text{m/s}$. We consider 49.5 $\text{m}$ long and 18 m depth domain for the test configuration. 
Figure \ref{true_fig} shows the S-wave velocity and P-wave velocity profiles for the assumed model. The soil layer (cyan color) is located approximately 7 m depth from the surface and the limestone layer (yellow color) is located from 8 m to 18 m depth. The void, blue rectangle in the domain, is located at the 15 m in the $x$ direction and 7 m in the $ z$ direction (depth). 

\begin{figure}
\centering
\includegraphics[width=1\linewidth]{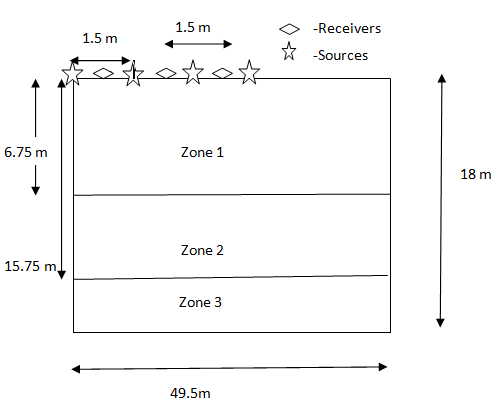}
\caption{The domain categorization as zone 1, zone 2, and zone 3. Receivers and sources are on the ground surface. }
\label{zone}
\end{figure}

\begin{figure}
\centering
\includegraphics[width=1\linewidth]{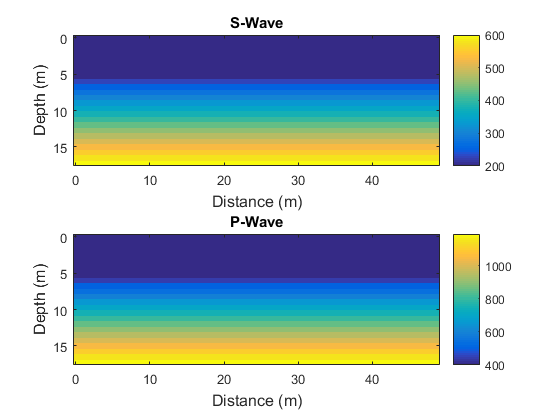}
\caption{Velocity profiles of the initial model}
\label{ini_fig}
\end{figure}

The finite difference code developed by Ref.~\cite{tran2012site} was used to generate a synthetic waveform data set. The code was modified to be used with the difference cell size method. The synthetic waveform data were
recorded from 32 receivers spaced every 1.5 m from station 0.75 m to 49.5 m. 33 shots were used at 1.5 m spacing starting from 0 to
36 m on the ground surface.  Fig.~\ref{zone} shows the receiver locations and source positions. The waveform data obtained with the finite difference code is used for inversion. For the data inversion, an initial model is generated with S-wave velocity increasing with depth (from 200 $\text{m/s}$ at the surface to 600 $\text{m/s}$ at the bottom) and P-wave velocity was generated from the S-wave velocity using Eq.~\ref{pwave}. Figure~\ref{ini_fig} shows the initial model, which was used for the inversion. Step size in the regular grid is $dx = dz =0.3$ in both $x$ and $z$ directions. Widths of the three regions for the difference cell size methods are 6.75 m, 9  m, and 2.25 m for zone 1, zone 2, and zone 3, respectively (see Fig.~\ref{zone}). The corresponding step sizes are $dx, 2dx,$ and $3dx$ for zone 1, zone 2, and zone 3, respectively. 

With the initial model, four inversions are performed for the data sets at four frequency ranges at central frequencies of 10, 15, and $20$ Hz. The first inversion at a central frequency of $10$ Hz started with the initial model. The other inversions at the central frequencies 15 and $20$ Hz were performed by using the inversion results at the lower central frequency as the initial model. During the inversion, S-wave and P-wave velocities were updated using Eq.~\ref{update}.

\begin{figure}
\centering
\includegraphics[width=0.8\linewidth]{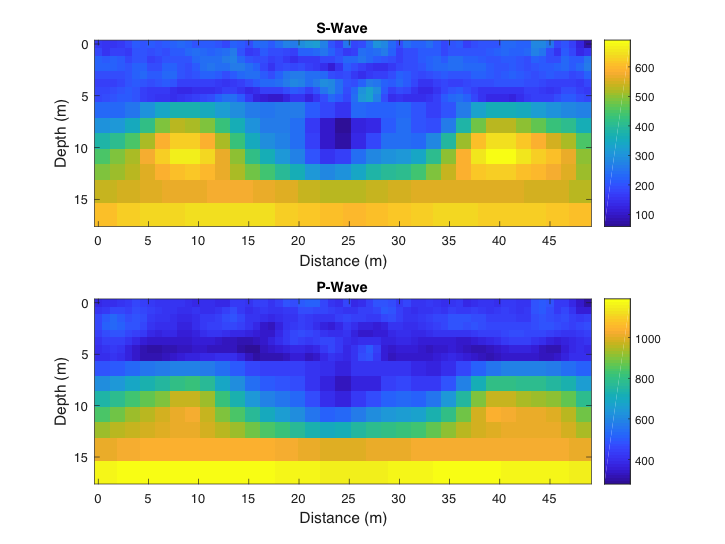}(a)
\includegraphics[width=0.8\linewidth]{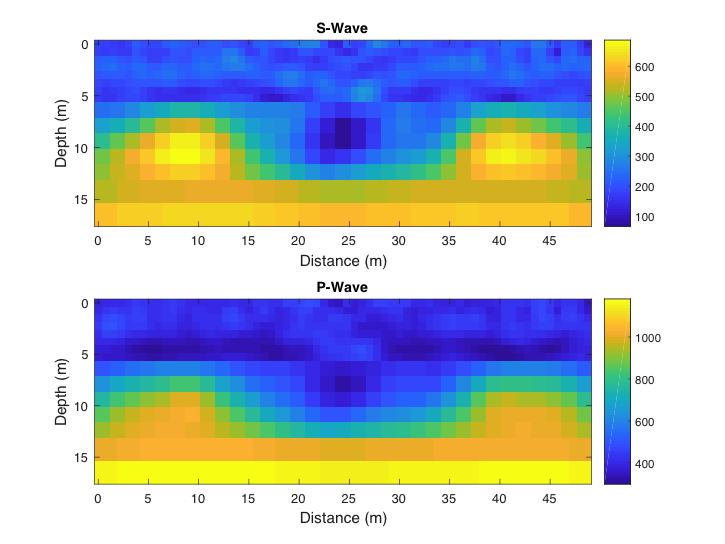}(b)
\caption{The inversion results for S-wave and P-wave velocities at the central frequency (a) $10$ Hz and (b) $15$ Hz.}
\label{fig_10_blo}
\end{figure}

The inversion results with the central frequency $10$ Hz and $15$ Hz are shown in Fig.~\ref{fig_10_blo}(a) and (b), respectively. At $10$ Hz, the void and the two layers can be clearly characterized by the S-wave velocity profile. Two layers can also be characterized from the P-wave velocity profile, but the void cannot be seen clearly from the P-wave velocity profile.  From the inversion results at $15 Hz$,  the void can be characterized by both S-wave and P-wave velocity profiles. Notice that inversion results at $10$ Hz were used as the initial model for the inversion at $15$ Hz.  

\begin{figure}
\centering
\includegraphics[width=1\linewidth]{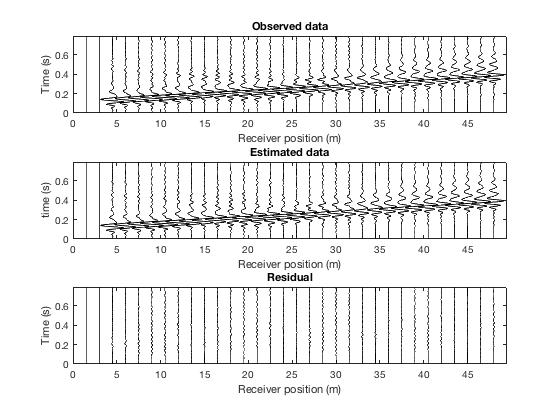}
\caption{The observed and estimated data for the inversion at the central frequency 20 Hz using the different cell size approach}
\label{fig_res}
\end{figure}

Convergence of the iteration method was tested by using the residual between the estimated and observed data. In all inversions, the convergence occurred at 20 iterations.  Figure \ref{fig_res} shows the estimated and observed waveforms at receiver positions for the inversion at the central frequency at 20 Hz. The residuals at the receivers are very small due the similar waveform of observed and estimated data. 

\begin{figure}
\centering
\includegraphics[width=1\linewidth]{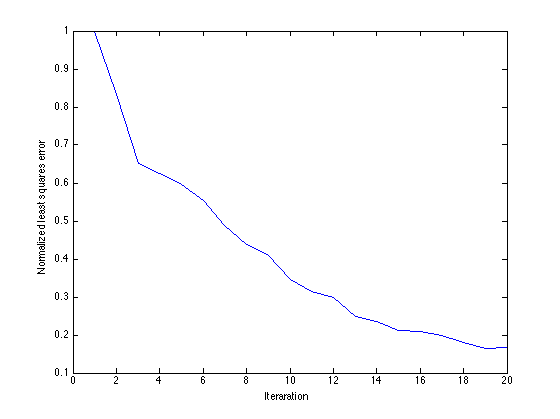}
\caption{Least square error as a function of integration number at simulations}
\label{fig_er}
\end{figure}

The normalized least squares error for 20 iterations are shown in Fig.~\ref{fig_er}. One can see the 0.8 order reduction in the error from the $1^{\text{st}}$ iteration to the $20^{\text{th}}$ iteration. After the $20^{\text{th}}$ iteration, the error reached a plateau and results are converged at the $20^{\text{th}}$ iteration. 

\subsubsection{Computational Efficiency of the Different Cell Size Method}

\begin{figure}
\centering
\includegraphics[width=0.8\linewidth]{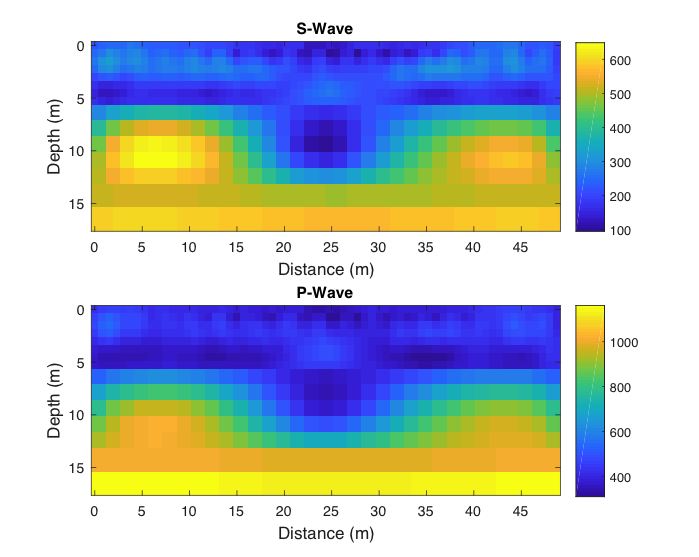}(a)
\includegraphics[width=0.8\linewidth]{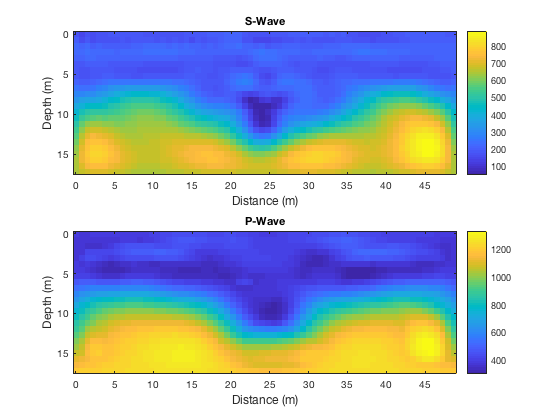}(b)
\caption{The inversion results for S-wave and P-wave velocities at the central frequency 20Hz using (a) the different cell size approach and (b) the regular grid method}
\label{fig_compare}
\end{figure}

The results of the different cell size method are compared with the regular cell size method used in Ref.~\cite{tran2012site}. The comparison here is done only to see the accuracy and computational efficiency of the different cell size method. The model updates at $20$ Hz from the regular grid method and different cell size method are shown in Fig.~\ref{fig_compare}. By comparing both models with the true model, one can accurately identify the void and the soil layers. Moreover, both methods are able to characterize the location, shape, and the S-wave velocity of the void. However, one should notice that the Hessian approximation matrix calculation with regular cell size method took about 3 hours on a Mac computer with a 2.6 GHz processor, while the different cell size method took only about 2.5 hours.  For the synthetic model that we considered here, the domain was discretized to have $24 \times 66$ =1584  cells. Thus the number of parameters in the model is 1584. Since the number of sources multiplied by the number of receivers is 1506, the size of the Jacobian matrix is $1506 \times 1584$ and the size of the $H_a = J^T J$ is 1584 $\times$ 1584. When we use our different cell size approach to combine the bottom cells of the domain, the number of cells in the domain was reduced to 814, so the number of parameters of the problem was reduced to 814. With this different cell size method, the size of the $H_a = J^T J$ is $814 \times 814$, which is less expensive to calculate. One can see that the size of the $H_a$ has reduced to approximately to 1/4 of the original $H_a$ matrix. The reduction of the size of $H_a$ depends on the discretization ratio of zone 3 and zone 2 and the decomposition of the domain. For real experiments, usually, the inversion problems are large scale problems.  In the 3-D problems, the size of the data sets and the number of parameters of the model are large. The different cell size method is competitive even with high resolutions.  For example, consider the case with $100 \times 150 = 15000$ cells, widths of the zone 1, zone 2, zone 3 are 46 m, 24 m, and 30 m, and discretization ratios for zone 1, zone 2, zone 3 are 1,2, and 3. Then the size of the new $H_a$ matrix is $8300 \times 8300$. The size of new $H_a$ has reduced approximately to 1/4 of the original $H_a$ matrix. 
The results of the difference cell size method are compared with the regular cell size method used in Ref.~\cite{tran2012site}. The comparison here is done only to see the accuracy and the computational efficiency of the different cell size method. The model updates at $20$ Hz from the regular grid method and different cell size method are shown in Fig.~\ref{fig_compare}. By comparing both models with the true model, one can accurately identify the void and the soil layers. Moreover, both methods are able to characterize the location, shape, and the S-wave velocity of the void. However, one should notice that, the Hessian approximation matrix calculation with regular cell size method took about 3 hours on a Mac computer with a 2.6 GHz processor, while the difference cell size method took only about 2.5 hours.  For the synthetic model that we considered here, the domain was discretized to have $24 \times 66$ =1584  cells. Thus the number of parameters in the model is 1584. Since the number of sources multiplied by the number of receivers is 1506, the size of the Jacobian matrix is $1506 \times 1584$ and the size of the $H_a = J^T J$ is 1584 $\times$ 1584. When we use our different cell size approach to combine the bottom cells of the domain, the number of cells in the domain was reduced to 814, so the number of parameters of the problem was reduced to 814. With this different cell size method, the size of the $H_a = J^T J$ is $814 \times 814$, which is less expensive to calculate.  One can see that the size of the $H_a$ has reduced to approximately to 1/4 of the original $H_a$ matrix. The reduction of the size of $H_a$ depends on the discretization ratio of zone 3 and zone 2 and the decomposition of the domain. For real experiments, usually, the inversion problems are large scale problems.  In the 3-D problems, the size of the data sets and the number of parameters of the model are large. The difference cell size method is competitive even with high resolutions. For example, consider the case with $100 \times 150 = 15000$ cells, widths of the zone 1, zone 2, zone 3 are 46 m, 24 m, and 30 m, and discretization ratios for zone 1, zone 2, zone 3 are 1,2, and 3. Then the size of the new $H_a$ matrix is $8300 \times 8300$. The size of new $H_a$ has reduced approximately to 1/4 of the original $H_a$ matrix. 
Thus the difference cell size method is more efficient than the regular grid method and has a good potential for 3-D wave inversion and large scale problems.  cell size method is more efficient than the regular grid method and has a good potential for 3-D wave inversion and large scale problems. 

\begin{figure}
\centering
\includegraphics[width=0.8\linewidth]{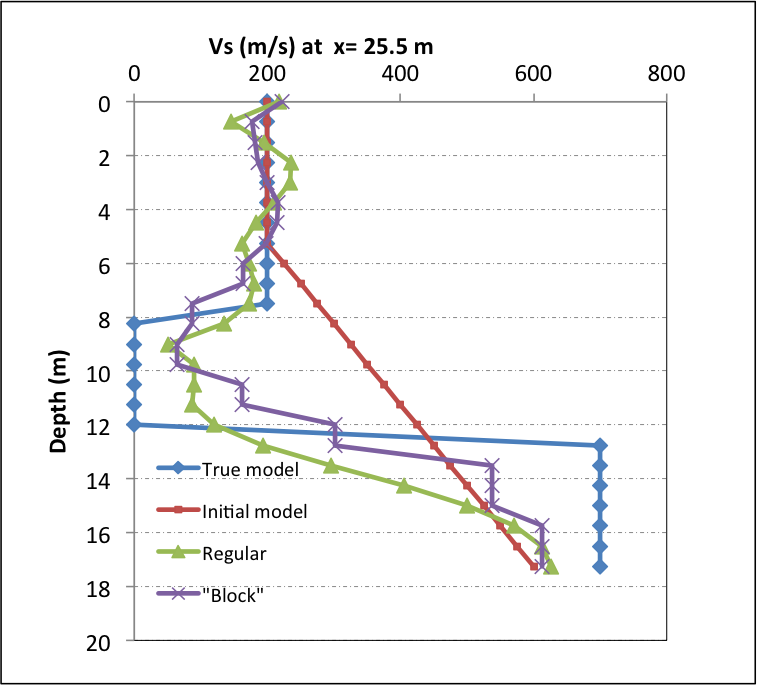}(a)
\includegraphics[width=0.8\linewidth]{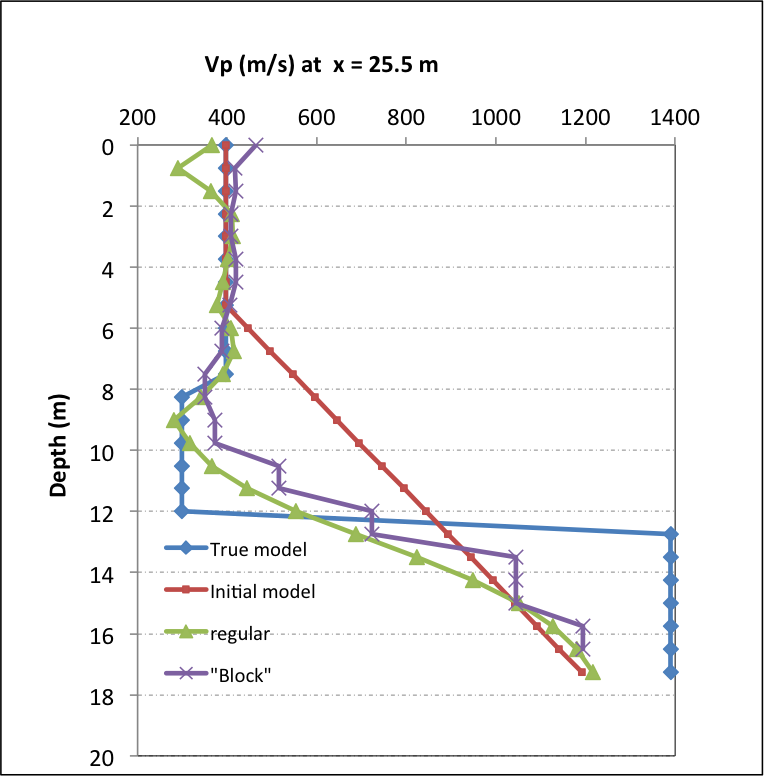}(b)

\caption{The variation of the velocity profile at $x = 25.5$ m (a) S-Wave velocity profile (b) P-Wave velocity profile}
\label{fig_vel}
\end{figure}

Figure \ref{fig_vel}(a) and (b) show the inverted 2-D profiles of S-wave velocity and  P-wave velocity variations along the depth, respectively. Blue color curve represents velocity variation with depth in the pre-assumed true model. The observed velocity variation from inversion using the regular method and the different cell method are shown in green color and purple. 
The first layer appears from 0 to 8 m depth. The second layer appears from 12 m to 18 m depth. The void is located from 8 m to 12 m depth. The velocity variations from both regular and difference cell size method closely follow the same variations as true model.  One can see that two layers, including the void, are clearly characterized by both velocity profiles.

\section{Discussion and Conclusions}

In this work, full seismic waveform inversion method using the Gauss-Newton method was utilized for detection of embedded sinkholes. The forward problem for simulating seismic wave fields was solved using the velocity-stress staggered-grid finite difference method. A model update of the inversion method was performed with the Gauss-Newton method with the difference cell size method. 

One of the major disadvantages of the  Gauss-Newton model updating is the large amount of computational and memory required to calculate the Hessian approximation matrix and the Jacobian matrix. To overcome the computation and memory requirements, we used difference cell size approach for storing the Jacobian matrix.  The values of the Jacobian, which were obtained using partial derivatives of seismograms with respect to the parameters, at the bottom part of the cells in the domain take smaller values compared with the values of Jacobian at the upper part of the domain. In this approach, the domain is decomposed into three zones according to the values of the Jacobian at the cells. The cells in the bottom of the zones are combined to create larger cells. The Jacobian matrix is then recalculated appropriately. 

The results are validated for a synthetic model with an embedded air-filled void. The synthetic model, which consists of two layers of an earth model, is investigated. The inversion at three frequency ranges with central frequencies of 10, 15, and 20 Hz were performed. The void can be characterized from both S-Wave and P-Wave velocities.
The computational requirements for both the difference cell size method and the regular cell size method are compared. The difference cell method is able to compute the Hessian with less computational time than required for regular method. In conclusion, the developed approach is well suitable for 3D full wave inversion with other Geo-technical conditions.

\begin{acknowledgments}
We thank Prof. Erik Bollt at Clarkson University, Potsdam and Prof. Khiem Tran at University of Florida, Florida for their valuable discussions, guidance, and comments.
\end{acknowledgments}


\section*{References}
\bibliography{intro,introbib}
\end{document}